\documentclass[12pt]{article}
\usepackage{graphicx}

\usepackage{epsfig}
\usepackage{amsmath}
\usepackage{amssymb}
\usepackage{amsthm}
\usepackage{latexsym}
\usepackage{cite}

\newcommand{\comment}[1]{}
\newtheorem{theorem}{Theorem}

\newtheorem{proposition}{Proposition}[section]

\begin{document}

\title{\LARGE
{\bf Towards Conformal Invariance \\ and a Geometric Representation \\ of the 2D Ising Magnetization Field}
\thanks{Based on a talk by the author on joint work with C.M.~Newman and work in progress with C.~Garban
and C.M.~Newman given at the workshop Inhomogeneous Random Systems
2010 (Paris).} }

\author{
{\bf Federico Camia}
\thanks{Research supported in part by NWO Vidi grant 639.032.916.}\,
\thanks{E-mail: f.camia@vu.nl}\\
{\small \sl Department of Mathematics, Vrije Universiteit Amsterdam}\\
{\small and}\\
{\small \sl NYU Abu Dhabi}
}

\date{}

\maketitle

\begin{abstract}
We study the continuum scaling limit of the critical Ising
magnetization in two dimensions. We prove the existence of
subsequential limits, discuss connections with the scaling limit of
critical FK clusters, and describe work in progress of the author
with C. Garban and C.M. Newman.
%
%
\end{abstract}



\section{Synopsis} \label{sec-syn}

The Ising model in $d=2$ dimensions is perhaps the most studied
statistical mechanical model and has a special place in the theory
of critical phenomena since the groundbreaking work of Onsager~\cite{onsager}.
Its scaling limit at or near the critical point is recognized to
give rise to Euclidean (quantum) field theories. In particular,
at the critical point, the lattice magnetization field should
converge, in the scaling limit, to a Euclidean random field $\Phi^0$
corresponding to the simplest reflection-positive conformal field
theory~\cite{bpz2,cardy08}.
As such, there have been a variety of representations in terms of
free fermion fields~\cite{sml} and explicit formulas for correlation
functions (see, e.g., \cite{mw-book,palmer} and references therein).

In~\cite{cn-pnas}, C.M. Newman and the present author introduced a
representation of $\Phi^0$ in terms of random geometric objects
associated with Schramm-Loewner Evolutions (SLEs)~\cite{schramm}
(see also~\cite{cardy3,kn,lawler2,werner4}) and Conformal Loop
Ensembles (CLEs)~\cite{sheffield,werner3,sw1,sw2}---namely, a gas
(or random process) of continuum loops and associated clusters and
(renormalized) area measures.

The purpose of the present paper is twofold, as we now explain.
First of all, we provide a detailed proof of the existence of
subsequential limits of the lattice magnetization field as a square
integrable random variable and a random generalized function
(Theorem~\ref{thm:subseqential-limits}) following the ideas
presented in~\cite{cn-pnas}. We also introduce a cutoff field whose
scaling limit admits a geometric representation in terms of rescaled
counting measures associated to critical FK clusters, and show that
it converges to the magnetization field as the cutoff is sent to
zero (Theorem~\ref{thm:cutoff-removal}).

Secondly, we describe work in progress~\cite{cgn} of the author with
C. Garban and C.M. Newman aimed at establishing uniqueness of the
scaling limit of the lattice magnetization and conformal covariance
properties for the limiting magnetization field. We also explain how
the existence and conformal covariance properties of the
magnetization field should imply the convergence, in the scaling
limit, of a version of the model with a vanishing (in the limit)
external magnetic field to a field theory with exponential decay of
correlations, and how they can be used to determine the free energy
density of the model up to a constant
(equation~(\ref{eq:free-energy})).

\section{The Magnetization and Some Results} \label{sec-intro}

We consider the standard Ising model on the square lattice ${\mathbb
Z}^2$ with (formal) \emph{Hamiltonian}
\begin{equation} \label{eq:hamiltonian}
{\bf H} = -\sum_{\{x,y\}} S_x S_y - H \sum_x S_x \, ,
\end{equation}
where the first sum is over nearest-neighbor pairs in ${\mathbb
Z}^2$,
the spin variables $S_x, S_y$ are $(\pm 1)$-valued and the external
field $H$ is in $\mathbb R$. For a bounded $\Lambda \subset {\mathbb
Z}^2$, the \emph{Gibbs distribution} is given by
$\frac{1}{Z_{\Lambda}} \, e^{-\beta \, {\bf H}_{\Lambda}}$, where
${\bf H}_{\Lambda}$ is the Hamiltonian~(\ref{eq:hamiltonian}) with
sums restricted to sites in $\Lambda$, $\beta \geq 0$ is the
\emph{inverse temperature}, and the \emph{partition function}
$Z_{\Lambda}$ is the appropriate normalization needed to obtain a
probability distribution.

We are mostly interested in the model with zero (or vanishing)
external field, and at the critical inverse temperature, $\beta_c =
\frac{1}{2} \, \log{(1+\sqrt{2})}$. For all $\beta\leq\beta_c$, the
model has a unique \emph{infinite-volume Gibbs distribution} for any
value of the external field $H$, obtained as a weak limit of the
Gibbs distribution for bounded $\Lambda$ by letting $\Lambda
\uparrow {\mathbb Z}^2$. For any value of $\beta\leq\beta_c$ and of
$H$, expectation with respect to the unique infinite-volume Gibbs
distribution will be denoted by $\langle \cdot \rangle_{\beta,H}$.
At the \emph{critical point}, that is when $\beta=\beta_c$ and
$H=0$, expectation will be denoted by $\langle
\cdot \rangle_c$. 
By translation invariance, the \emph{two-point correlation} $\langle
S_x S_y \rangle_{\beta,H}$ is a function only of $y-x$, which at the
critical point we denote by $\tau_c(y-x)$.


We want to study the random field associated with the spins on the
rescaled lattice $a \, {\mathbb Z}^2$ in the scaling limit $a \to
0$. More precisely, for functions $f$ of bounded support on
${\mathbb R}^2$, we define for the critical model
\begin{equation} \label{eq:lat-field1}
\Phi^a(f) \equiv \int_{{\mathbb R}^2} f(z) \Phi^a(z) dz \equiv
\int_{{\mathbb R}^2} f(z) [ \Theta_a \sum_{x \in {\mathbb Z}^2} S_x
\delta(z-ax) ] dz = \Theta_a \sum_{z \in a \, {\mathbb Z}^2} f(z)
S_{z/a} \, ,
\end{equation}
with scale factor
\begin{equation} \label{eq:Theta-first}
\Theta_a^{-1} \equiv \sqrt{ \sum_{z,w \in \Lambda_{1,a}} \langle
S_{z/a} S_{w/a} \rangle_c } = \sqrt{ \sum_{x,y \in \Lambda_{1/a}}
\tau_c(y-x) } \, ,
\end{equation}
where $\Lambda_{L,a} \equiv [0,L]^2 \cap a \, {\mathbb Z}^2$ and
$\Lambda_L \equiv \Lambda_{L,1} = [0,L]^2 \cap {\mathbb Z}^2$.

The block magnetization, $M^a \equiv \Phi^a({\bf 1}_{[0,1]^2})$,
where $\bf 1$ denotes the indicator function, is a rescaled sum of
identically distributed, \emph{dependent} random variables. In the
high temperature case, $\beta < \beta_c$, and with zero external
field, $H=0$, the dependence is sufficiently weak for the block
magnetization to converge, as $a \to 0$, to a mean-zero, Gaussian
random variable (see, e.g.,~\cite{newman80} and references therein).
In that case, the appropriate scaling factor $\Theta_a$ is of order
$a$, and the field converges to Gaussian white noise as $a \to 0$
(see, e.g.,~\cite{newman80}). In the critical case, however,
correlations are much stronger and extend to all length scales, so
that one does not expect a Gaussian limit. A proof of this will be
presented elsewhere \cite{cgn}; in this paper we are concerned with
the existence of subsequential limits for the lattice magnetization
field, and their geometric representation in terms of area measures
of critical FK clusters.

\begin{figure}[!ht]
\begin{center}
\includegraphics[width=9cm]{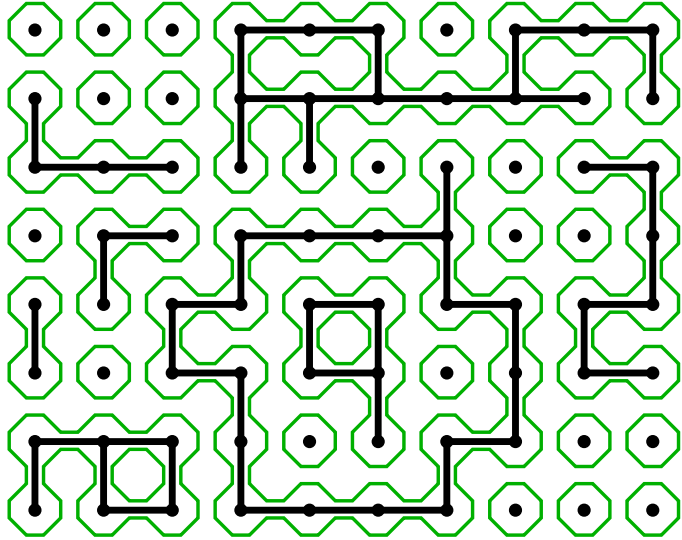}
\caption{Example of an FK bond configuration in a rectangular region.
Black dots represent sites of ${\mathbb Z}^2$, black horizontal and vertical
edges represent FK bonds. The FK clusters are highlighted by lighter (green)
loops on the medial lattice.}
\label{fig:loops}
\end{center}
\end{figure}

The FK representation of the Ising model with zero external field,
$H=0$, is based on the $q=2$ random-cluster measure $P_p$
(see~\cite{grimmett-rcm-book} for more on the random-cluster model
and its connection to the Ising model). A spin configuration
distributed according to the unique infinite-volume Gibbs
distribution with $H=0$ and inverse temperature $\beta\leq\beta_c$
can be obtained in the following way. Take a random-cluster (FK)
bond configuration on the square lattice distributed according to
$P_p$ with
$p=p(\beta)= 1 - e^{-2\beta}$, and let $\{ {\cal C}_i \}$ 
denote the corresponding collection of FK clusters, where a cluster
is a maximal set of sites of the square lattice connected via bonds
of the FK bond configuration (see Figure~\ref{fig:loops}).
One may regard the index $i$ as taking values in the natural numbers,
but it's better to think of it as a dummy countable index without
any prescribed ordering, like one has for a Poisson point process.
Let
$\{ \eta_i \}$ 
be ($\pm 1$)-valued, i.i.d., symmetric random  variables, and
assign $S_x=\eta_i$ for all $x \in {\cal C}_i$; then the collection
$\{ S_x \}_{x \in {\mathbb Z}^2}$ of spin variables is distributed
according to the unique infinite volume Gibbs distribution with $H=0$
and inverse temperature $\beta$. When $\beta=\beta_c$, we will use
the notation $P_c \equiv P_{p(\beta_c)}$, and $E_c$ for expectation
with respect to $P_c$.

A useful property of the FK representation is that, when $H=0$, the
Ising two-point function can be written as
\begin{equation} \nonumber
\langle S_x S_y \rangle_{\beta,0} = P_{p(\beta)}(x \text{ and } y
\text{ belong to the same FK cluster } {\cal C}_i) \, .
\end{equation}
As an immediate consequence, we have
\begin{equation} \label{eq:Theta}
\Theta^{-2}_a = \sum_{x,y \in \Lambda_{1/a}} \tau_c(y-x) = \sum_{x,y
\in \Lambda_{1/a}} E_c \left[ \sum_i {\bf 1}_{x \in {\cal C}_i} {\bf
1}_{y \in {\cal C}_i} \right] = E_c \left[ \sum_i | \hat{\cal C}^a_i
|^2 \right] \, ,
\end{equation}
where $\hat{\cal C}^a_i$ is the restriction of the rescaled cluster
${\cal C}^a_i=a \, {\cal C}_i$ in $a \, {\mathbb Z}^2$ to $[0,1]^2$,
and $|\hat{\cal C}^a_i|$ is the number of ($a \, {\mathbb
Z}^2$)-sites in $\hat{\cal C}^a_i$. (Note that $\hat{\cal C}^a_i$
need not be connected.) Using the FK representation, we can
write~(\ref{eq:lat-field1}) as
\begin{equation} \label{eq:lattice-field}
\Phi^a(f) \stackrel{dist.}{=} \sum_i \eta_i \mu^a_i(f) \, ,
\end{equation}
where $\mu^a_i \equiv \Theta_a \sum_{x \in {\cal C}_i}\delta(z-ax) $
and the $\eta_i$'s, as before, are $(\pm 1)$-valued, symmetric
random variables independent of each other and everything else. We
can now easily see that $\Theta_a$ was chosen so that the second
moment of the block magnetization $M^a$, defined earlier, is exactly
one:
\begin{equation} \label{eq:second-moment}
\left\langle \left( M^a \right)^2 \right\rangle_c = \left\langle
\left[ \Phi^a({\bf 1}_{[0,1]^2}) \right]^2 \right\rangle_c =
E_c\left[\sum_i \left(\mu^a_i({\bf 1}_{[0,1]^2})\right)^2\right] =
\Theta^2_a E_c \left[ \sum_i |\hat{\cal C}^a_i|^2 \right] = 1 \, .
\end{equation}


We can associate in a unique way to each rescaled counting measure
$\mu^a_i$ the interface $\gamma^a_i$ in the medial lattice between
the corresponding (rescaled) FK cluster ${\cal C}^a_i$ and the
surrounding FK clusters. Since all FK clusters are almost surely
finite at the critical point ($\beta=\beta_c, H=0$), such interfaces
form closed curves, or loops, which separate the corresponding
clusters ${\cal C}^a_i$ from infinity (see Fig.~\ref{fig:loops}).
There are two types of loops: (1) those with sites of $a \, {\mathbb
Z}^2$ immediately on their inside and (2) those with sites of $a \,
{\mathbb Z}^2$ immediately on their outside.
We denote by $\{ \gamma^a_i \}$ the (random) collection of all loops
of the first type associated with the FK clusters $\{ {\cal C}^a_i
\}$. Each realization of $\{ \gamma^a_i \}$ can be seen as an
element in a space of collections of loops with the
Aizenman-Burchard metric~\cite{ab}. (The latter is the induced
Hausdorff metric on collections of curves associated to the metric
on curves given by the infimum over monotone reparametrizations of
the supremum norm.)
It follows from~\cite{ab} and the RSW-type bounds of~\cite{dhn} (see
Section~5.3 there) that, as $a \to 0$, $\{ \gamma^a_i \}$ has
subsequential limits in distribution to random collections of loops
in the Aizenman-Burchard metric. In the scaling limit, one gets
collections of nested loops that can touch (themselves and each
other), but never cross.

In order to study the magnetization field, we introduce some more notation.
Let $(C_0({\mathbb R}^2),||\cdot||_{\infty})$ denote the space
of continuous functions on ${\mathbb R}^2$ with compact support,
endowed with the metric of uniform convergence.
Let $({\cal P}_2,W_2)$ denote the space of probability distributions
on $\mathbb R$ (with the Borel $\sigma$-algebra) with finite second
moment, endowed with the Wasserstein (or minimal $L_2$) metric
\begin{equation} \label{eq:W-distance}
W_2(P,Q) \equiv \left( \inf E \left[ |X-Y|^2 \right] \right)^{1/2} \, ,
\end{equation}
where $X$ and $Y$ are coupled random variables with respective
distributions $P$ and $Q$, $E$ denotes expectation with respect
to the coupling, and the infimum is taken over all such couplings
(see, e.g., \cite{ruschendorff01} and references therein).
Convergence in the Wasserstein metric $W_2$ is equivalent to
convergence in distribution plus convergence of the second
moment. For brevity, we will write $C_0({\mathbb R}^2)$ and
${\cal P}_2$, instead of $(C_0({\mathbb R}^2),||\cdot||_{\infty})$
and $({\cal P}_2,W_2)$, unless we wish to emphasize the role
of the metrics.

We further denote by ${\cal D}$ the space of infinitely
differentiable functions on ${\mathbb R}^2$ with compact support,
equipped with the topology of uniform convergence of all
derivatives, and by ${\cal D}'$ its topological dual, i.e., the
space of all generalized functions.

The next theorem shows that the lattice magnetization field has
subsequential scaling limits in terms of continuous functionals,
in a distributional sense using the Wasserstein metric $W_2$,
and in the sense of generalized functions by an application of
the Bochner-Minlos theorem. (We remark that the last statement of
Theorem~\ref{thm:subseqential-limits} is not optimal in the sense
that similar conclusions should apply to a larger class of functions
than $\cal D$.)

\begin{theorem} \label{thm:subseqential-limits}
For any sequence $a_n \to 0$, there exists a subsequence $a_{n_k} \to 0$
such that, for all $f \in C_0({\mathbb R}^2)$, the distribution $P^k_f$
of $\Phi^{a_{n_k}}(f)$ converges in the Wasserstein metric~(\ref{eq:W-distance}),
as $k \to \infty$, to a limit $P^0_f \in {\cal P}_2$ such that the map
$P^0_{\cdot}: (C_0({\mathbb R}^2),||\cdot||_{\infty}) \longrightarrow ({\cal P}_2,W_2)$
is continuous.
Furthermore, for every subsequential limit $P^0_{\cdot}$,
there exists a random generalized function $\Phi^0 \in {\cal D}'$
with characteristic function $\chi(f) \equiv \int e^{ix} dP^0_f(x)$.
\end{theorem}

Theorem~\ref{thm:subseqential-limits} represents the starting point
of a joint project with C. Garban and C.M. Newman aimed at
establishing uniqueness of the scaling limit of the lattice
magnetization field and its conformal covariant properties. One not
only expects a unique scaling limit for the lattice magnetization
field, but based on the representation~(\ref{eq:lattice-field}), one
would like to write the limiting field $\Phi^0$ as
\begin{equation} \label{eq:geom-rep}
``\Phi^0(f) = \sum_{j} \eta_j \mu^0_j(f)"
\end{equation}
where the $\mu^0_j(f)$'s are the putative scaling limits of the
$\mu^a_i(f)$'s that appear in~(\ref{eq:lattice-field}). Indeed, in
the scaling limit, one should obtain a collection $\{ \mu^0_j \}$ of
mutually orthogonal, finite measures supported on the scaling limit
of the critical FK clusters. However, due to scale invariance, $\{
\mu^0_j \}$ should contain (countably) infinitely many elements, and
the scaling covariance expected for the $\mu^0_j$'s suggests that
the collection $\{\mu^0_j(f)\}$ is in general not absolutely
summable. What meaning, if any, can we then attribute to the sum
in~(\ref{eq:geom-rep})?

To help answer that question,
we introduce the $\varepsilon$-cutoff lattice magnetization field
\begin{equation} \label{eq:cutoff-lattice-field}
\Phi^a_{\varepsilon}(f) \equiv \sum_{i:
diam({\gamma}^a_i)>\varepsilon} \eta_i \mu^a_i(f) \, ,
\end{equation}
where the elements of the collection $\{ \mu^a_i \}$ of all rescaled
(random) measures that are involved
in~(\ref{eq:cutoff-lattice-field}) are those associated to rescaled
FK clusters ${\cal C}^a_i$ that intersect the support of $f$ and
whose corresponding loops $\gamma^a_i$ have diameter $>\varepsilon$.

Once again, one would like to write the scaling limit of the cutoff
field as ``$\Phi^0_{\varepsilon}(f) = \sum_{j:
diam({\gamma}^0_j)>\varepsilon} \eta_j \mu^0_j(f)$". In this case
however, the sum would be unambiguous because it would contain only
a \emph{finite} number of terms. A proof of the latter fact follows
from Prop.~\ref{prop:crossprob} in Section~\ref{sec:beyond}.
Combined with~(\ref{eq:second-moment}) and
Prop.~\ref{prop:small-clusters} in Section~\ref{sec:proofs},
Prop.~\ref{prop:crossprob} implies that the collection of
$\mu^a_i(f)$'s corresponding to macroscopic FK clusters has
nontrivial subsequential scaling limits. Indeed, it is clear from
equation~(\ref{eq:second-moment}) that no $\mu^a_i({\bf
1}_{[0,1]^2})$ can diverge as $a \to \infty$. In addition,
Prop.~\ref{prop:small-clusters} says that ``small'' FK clusters do
not contribute to the magnetization in the scaling limit and thus,
by Prop.~\ref{prop:crossprob}, the number of FK clusters which
contribute significantly to $M^a=\Phi^a({\bf 1}_{[0,1]^2})$ remains
\emph{bounded} as $a \to 0$.
Since $\left\langle (M^a)^2 \right\rangle_c = 1$ for all $a$,
this implies that not all
$\mu^a_i({\bf 1}_{[0,1]^2})$'s can converge to 0 as $a \to 0$.
Prop.~\ref{prop:tau} ensures that the same conclusions hold not only
for the collection of $\mu^a_i(f)$'s with $f={\bf 1}_{[0,1]^2}$, but
for other functions as well.

The result below shows that, in the scaling limit, one recovers
the ``full'' magnetization field from the cutoff one by letting the
cutoff go to zero.

\begin{theorem} \label{thm:cutoff-removal}
For any sequence $a_n \to 0$, there exists a subsequence $a_{n_k}
\to 0$ such that, for all $f \in C_0({\mathbb R}^2)$ and all $m \in
{\mathbb N}$, the distributions of $\Phi^{a_{n_k}}(f)$ and
$\Phi^{a_{n_k}}_{1/m}(f)$ converge in the Wasserstein
metric~(\ref{eq:W-distance}) as $k \to \infty$. Moreover, if $P^0_f$
and $P^0_{f,m}$ denote the respective limits, $P^0_{f,m}$ converges
to $P^0_f$ in the Wasserstein metric~(\ref{eq:W-distance}) as $m \to
\infty$.
\end{theorem}
In view of Theorem~\ref{thm:cutoff-removal}, one can interpret the
sum in equation~(\ref{eq:geom-rep}) as a shorthand for the limit of
the cutoff field as the cutoff is removed.
Combined with the fact that the collection of $\mu^a_i(f)$'s has
nontrivial subsequential scaling limits, as explained above,
Theorems~\ref{thm:subseqential-limits} and~\ref{thm:cutoff-removal}
partly establish the geometric representation proposed
in~\cite{cn-pnas}. In order to establish the existence of a unique
scaling limit for the collection of $\mu^a_i$'s as measures, and to
obtain their conformal covariance properties and those of the
limiting magnetization field $\Phi^0$, more work is needed. This is
discussed in the next section.

\section{Work in Progress: Uniqueness and Conformal Covariance}
The lattice magnetization field is expected to have a unique scaling limit
$\Phi^0$ with the property of transforming covariantly under conformal
transformations, i.e., if $\varphi$ is a conformal map,
\begin{equation} \label{eq:conf-cov}
\Phi^0(\varphi(z)) \stackrel{dist.}{=} |\varphi'(z)|^{-1/8}
\Phi^0(z) \, ,
\end{equation}
where $1/8$ is the Ising magnetization exponent. (With an abuse of notation,
we identify ${\mathbb R}^2$ and the complex plane $\mathbb C$.)

It is natural to attempt to prove such results using announced
results for FK percolation (see~\cite{smirnov-icm,smirnov-ising1})
which identify the scaling limit of the FK cluster boundaries (see
Fig.~\ref{fig:loops}) with SLE-type random fractal curves whose
distribution is invariant under conformal transformations. In order
to exploit such results, one can use techniques developed
in~\cite{gps,garban} to study the scaling limit of Bernoulli and
dynamical percolation in two dimensions. Roughly speaking, the idea
is to prove that the scaling limit of the ensemble $\{ \mu^a_i \}$
of rescaled counting measures associated to the FK clusters is a
measurable function of the collection 
of limiting (macroscopic) loops between FK clusters. 

To illustrate the idea, we take a small detour and discuss briefly
the scaling limit of Bernoulli percolation, focusing on site
percolation on the triangular lattice. The ``full" scaling limit of
percolation, comprising all interface loops separating macroscopic
clusters, was obtained by Camia and Newman in~\cite{cn04,cn06} and
shown to be a (nested) Conformal Loop Ensemble (CLE) in~\cite{cn08}.
In~\cite{cfn1,cfn2} Camia, Fontes and Newman proposed to construct
the near/off-critical scaling limit of percolation, with density of
open sites $p = 1/2 + \lambda a^{3/4}$ (where $\lambda \in
(-\infty,\infty)$ is a parameter, $a$ the lattice spacing, and $3/4$
the percolation correlation length exponent), from the critical one
``augmented'' by a ``Poissonian cloud'' of marks on the double
points of the limiting loops (i.e., where a loop touches itself or
where two different loops touch each other). Back on the lattice,
the marked points would correspond to ``pivotal'' sites that switch
state when the density of open sites is changed from $1/2$ to $p$,
causing a macroscopic change in connectivity. (The last sentence
should be interpreted in the context of the canonical coupling of
percolation models at different densities of open sites. In this
coupling, a percolation model with density $p$ of open sites is
obtained by assigning independent, uniform random variables $u_x \in
[0,1]$ to the sites $x$ of the lattice, and declaring open all sites
with $u_x<p$, and closed all other sites.) A key step in the
implementation of this idea is the construction of the intensity
measure of the Poisson process of marks. Since the points to be
marked are double points, it was argued in~\cite{cfn1,cfn2} that the
intensity measure should arise as the scaling limit of the
appropriately rescaled counting measure of
$\varepsilon$-macroscopically pivotal sites on the lattice with
spacing $a$, where an $\varepsilon$-macroscopically pivotal site $x$
has four neighbors which are the starting points of four alternating
paths, two made of (nearest-neighbor) open sites and two of closed
ones, reaching a distance $\varepsilon$ away from $x$.

The occurrence of an $\varepsilon$-macroscopically pivotal site $x$
in a percolation configuration is called a \emph{four-arm event}.
The scaling limit of the counting measure of
$\varepsilon$-macroscopically pivotal sites was obtained by Garban,
Pete and Schramm~\cite{gps} (see also~\cite{garban}) and used by the
same authors, in the spirit of the program proposed by Camia, Fontes
and Newman, to construct the near/off-critical scaling limit of
percolation. In particular, Garban, Pete and Schramm~\cite{gps}
consider the joint distribution of the collection of interface loops
and the (random) counting measure of $\varepsilon$-macroscopically
pivotal sites, $(\{ \gamma^a_i \}, \lambda^a_{\varepsilon})$, and
show that it converges to the law of some random variable $(\{
\gamma_j^0 \}, \lambda^0_{\varepsilon})$, where $\{ \gamma_j^0 \}$
is the collection of limiting loops and $\lambda^0_{\varepsilon}$ is
a random Borel measure. Moreover, they show that
$\lambda^0_{\varepsilon}$ is a measurable function of $\{ \gamma_j^0
\}$.

This last observation is in fact crucial, since the known uniqueness
of the scaling limit of the interface loops implies the uniqueness
of $\lambda^0_{\varepsilon}$. In addition, one can deduce how
$\lambda^0_{\varepsilon}$ changes under conformal transformations
from the knowledge of how $\{ \gamma_j^0 \}$ changes under those
same transformations. The latter can be deduced for the collection
$\{ \gamma_j^0 \}$ from the fact that it is a nested CLE whose loops
are SLE-type curves.

Heuristically, one can convince oneself that it is reasonable to
expect that $\lambda^0_{\varepsilon}$ be a measurable function of
$\{ \gamma_j^0 \}$ by noticing that knowing the macroscopic loops
should be sufficient to give a good estimate of the number of
macroscopically pivotal sites. For a discussion on how to turn this
observation into a proof, the reader is referred to Sect.~4.3
of~\cite{gps}, where complete proofs of the results mentioned in the
previous paragraph can also be found.

In Sect.~5 of~\cite{gps}, the authors discuss how to obtain similar
results for rescaled counting measures of other special sites. In
particular, they show how to obtain what they call the ``cluster''
or ``area'' measure, which counts the number of open sites contained
in clusters of diameter larger than some cutoff $\varepsilon>0$. The
occurrence of such a site $x$ corresponds to the event that there is
a path of (nearest-neighbor) open sites starting at $x$ and reaching
a distance $\varepsilon$ away from $x$. Such an event is called a
\emph{one-arm event}, and we will call $x$ a \emph{one-arm site}.
The proof in this case is in fact simpler because the event is
simpler, involving only one path.

At this point the reader should note that the area measures
$\mu^a_i$ introduced in the previous section in connection with the
magnetization field also count one-arm sites, with the only
difference that the relevant one-arm events are now in the context
of FK bond percolation.
FK percolation is more difficult to analyze than Bernoulli
percolation, due to the dependencies in the distribution of FK
configurations (as opposed to the product measure corresponding to
Bernoulli percolation). However, it seems that one can successfully
adapt the techniques of~\cite{gps,garban}, at least for the case of
one-arm sites which is relevant for the magnetization. As a
consequence, thanks to the results announced
in~\cite{smirnov-icm,smirnov-ising1}, one should obtain uniqueness
of the limiting ensemble $\{ \mu^0_j \}$ of area measures for the FK
clusters and of the magnetization field $\Phi^0$, as well as a proof
of~(\ref{eq:conf-cov}) and of the fact that, for any conformal map
$\varphi$, $\{|\varphi'(z)|^{-15/8}d\mu_j^0(\varphi(z))\}$ is
equidistributed with $\{d\mu_j^0(z)\}$. Because of the latter
property, we call the putative collection of measures $\{ \mu^0_j
\}$, obtained as the scaling limit of the collection of rescaled
counting measures $\{ \mu^a_i \}$, a \emph{Conformal Measure
Ensemble}.


\section{More Work in Progress: Free Energy Density and Tail Behavior}
The uniqueness and conformal covariance properties of $\Phi^0$  play
an important role in the analysis of the \emph{near-critical}
scaling limit (called \emph{off-critical} in the physics literature)
with a vanishing (in the limit) external field (at the critical
inverse temperature $\beta_c$). More precisely, consider an Ising
model on $a \, {\mathbb Z}^2$ with (formal)
Hamiltonian~(\ref{eq:hamiltonian}) and external field $H(a) = h
\beta_c^{-1} \Theta_a$ inside the square $[-L,L]^2$, and zero
outside it. We call $h$ the \emph{renormalized external field} and
note that the term
\begin{equation} \nonumber
-h \beta_c^{-1} \Theta_a \sum_{z \in a \, {\mathbb Z}^2 \cap [-L,L]^2} S_{z/a}
\end{equation}
in the Hamiltonian implies that the Gibbs distribution of this
particular Ising model is given by
\begin{equation} \nonumber
d\nu_{h,L}^a \equiv \frac{1}{Z_{h,L}^a} \exp\left(h \Theta_a \sum_{z
\in a \, {\mathbb Z}^2 \cap [-L,L]^2} S_{z/a} \right) d\nu^a =
\frac{1}{Z_{h,L}^a} \exp\left(h M^a_L \right) d\nu^a \, ,
\end{equation}
where $\nu^a$ is the Gibbs distribution corresponding to zero
external field, $Z_{h,L}^a$ is the appropriate normalization factor,
and $M^a_L$ denotes the block magnetization inside $[-L,L]^2$. As a
consequence, in the scaling limit ($a \to 0$) one would obtain a
distribution $\nu^0_{h,L}$ such that
\begin{equation} \nonumber
d\nu^0_{h,L} \equiv \frac{1}{Z^0_{h,L}} \exp\left(h \Phi^0({\bf
1}_{[-L,L]^2})\right) d\nu^0 \, ,
\end{equation}
where $Z^0_{h,L} \equiv \int \exp(h \Phi^0({\bf 1}_{[-L,L]^2})
d\nu^0$ and $\nu^0$ is the limiting distribution corresponding to
zero external field.

The question is now whether $\nu^0_{h,L}$ converges to some
$\nu^0_h$ as $L \to \infty$, and whether $\nu^0_h$ corresponds to
the physically correct near/off-critical scaling limit.
Heuristically, the correct normalization to obtain a nontrivial
near/off-critical scaling limit is such that the correlation length
$\xi$ remains bounded away from zero and infinity. Scaling theory
implies that $\xi \sim H^{-8/15}$ for small external field $H$. This
gives $H \sim a^{15/8}$, which coincides with the normalization
needed to obtain a nontrivial magnetization field (given by
$\Theta_a$), as can be seen from~(\ref{eq:Theta-first}) and the
asymptotic behavior of~$\tau_c$. With this in mind, we consider an
Ising model on $a \, {\mathbb Z}^2$ with an external field
$H=a^{15/8}$ inside $\Lambda_{L,a}$ and $0$ outside, for some large
$L$. Using the two-dimensional Ising critical exponent $\delta=15$
for the magnetization (i.e., $\langle S_0 \rangle_{\beta_c,H} \sim
H^{1/15}$ for small $H$, where $S_0$ denotes the spin at the
origin), and denoting by $\sum_x^L$ the sum over $x$ in $\Lambda_{L
/a}$, we can write the block magnetization in the unit square as
\begin{equation} \nonumber
\frac{\langle \Theta_a \sum_x^1 S_x \exp(a^{15/8} \sum_x^L S_x)
\rangle_c} {\langle \exp(a^{15/8} \sum_x^L S_x) \rangle_c}
\stackrel{L\gg1}{\sim} a^{15/8} a^{-2} \langle S_0
\rangle_{\beta_c,H=a^{15/8}} \sim a^{-1/8} (a^{15/8})^{1/15} = 1 \,
.
\end{equation}
Since the result is finite, this rough computation suggests a
positive answer to the previous question.

Indeed, using the convergence of the lattice magnetization field to
the continuum one and scaling properties of the critical FK
clusters, it appears possible to show~\cite{cgn} that, as $L \to
\infty$, $\nu^0_{h,L}$ has a unique weak limit, denoted by
$\nu^0_h$, and that $\nu^0_h$ represents the scaling limit of the
Ising model on $a \, {\mathbb Z}^2$ with external field $H(a) = h
\beta_c^{-1} \Theta_a$ on the whole plane.

The idea behind a proof of this makes use of the well-known ``ghost
spin'' representation of the Ising model with an external field, in
which an additional site with spin that agrees with the external
field is added and connected to all the sites of the square lattice.
The external field term in the Hamiltonian can then be written
(formally) as $-|H| \sum_x S_x S_g$, where the ghost spin $S_g$ is
equal to the sign of the external field $H$.
One can describe the Ising model with an external field using the FK
representation on the new graph comprising the square lattice and
the additional site carrying the ghost spin. Note however that the
density of FK bonds incident on the site carrying the ghost spin is
not given by $p(\beta) = 1 - e^{-2\beta}$, as for the other bonds,
but by $1 - e^{-2 \beta |H|}$.

The following key observation is an easy consequence of standard
properties of FK percolation. If a subset $\Lambda$ of the square
lattice is surrounded by a circuit $\Gamma$ of FK bonds that belong
to a cluster which also contains the site carrying the ghost spin,
the FK and spin configurations in $\Lambda$ are independent of the
FK and spin configurations outside the circuit $\Gamma$. The
RSW-type bounds proved in~\cite{dhn}, together with the FKG
inequality~\cite{fkg} and scaling properties of the FK clusters and
their area measures, imply that the probability to find such a
circuit $\Gamma$ surrounding any bounded subset $\Lambda$ is one.
This shows that the $\nu^0_{h,L}$-probability of any event that
depends only on the restriction of the spin configuration to a
finite subset $\Lambda$ of the square lattice has a limit as $L \to
\infty$. Consequently, the distribution $\nu^0_{h,L}$ has a weak
limit $\nu^0_h$ as $L \to \infty$.

It is interesting to note that the argument alluded to above also
shows that $\nu^0_h$ is locally absolutely continuous with respect
to the zero-field measure $\nu^0$. This is in contrast to the
situation in two-dimensional percolation, where the critical and
near-critical measures are mutually singular~\cite{nw}. It should be
noted, however, that the Ising analogue of that type of percolation
near-critical scaling limit is to set $H=0$ and let $\beta(a) \to
\beta_c$, rather than set $\beta = \beta_c$ and let $H(a) \to 0$.

One expects the near/off-critical field to be ``massive'' in the
sense that correlations under $\nu^0_h$ should decay exponentially
at large distances. To understand why this should be the case, it is
again useful to resort to the ghost spin representation discussed
earlier. Remember that the Ising two-point function can be expressed
in terms of connectivity properties of the FK clusters (see the
discussion about the FK representation preceding
equation~(\ref{eq:Theta})). Because of that, exponential decay of
correlations is equivalent to the statement that, if two sites of
the square lattice, $x$ and $y$, belong to the same FK cluster
${\cal C}_i$, the probability that ${\cal C}_i$ does not contain the
site carrying the ghost spin decays exponentially in the distance
between $x$ and $y$.
But the scaling law for the area measures, $d\mu_j^0(\alpha \, z)
\stackrel{dist.}{=} \alpha^{15/8} \, d\mu_j^0(z)$ for all
$\alpha>0$, suggests that a macroscopic FK cluster of diameter at
least $||x-y||=O(1)$ (that is, of order $a^{-1}$ in units of the
lattice spacing $a$) should contain at least $O(a^{-15/8})$ sites,
precisely enough to compensate for the small intensity of the
external field $H \sim a^{15/8}$, which determines the probability
of a cluster to contain the site carrying the ghost spin via the
density, $1 - e^{-2 \beta |H|}$, of FK bonds connected to that site.



The exponential decay of correlations can be used to show the
existence of the \emph{free energy density} $f(h)$ at the critical
(inverse) temperature, defined by
\begin{equation} \nonumber
f(h) \equiv -\beta_c^{-1}\lim_{L \to \infty} (2L)^{-2}
\log{\left(\int \exp(h \Phi^0({\bf 1}_{[-L,L]^2})) d\nu^0\right)} \,
,
\end{equation}
provided that the limit exists. (Because of symmetry, it suffices to
consider positive external fields, $h \geq 0$.) For the
nearest-neighbor lattice Ising model, following a standard argument
(see for instance~\cite{minlos99}, Lecture 8), one can show the
existence of the free energy by partitioning $[-L,L]^2$ into equal
squares of fixed size and writing the Hamiltonian as a sum of terms
of two types: those corresponding to the interactions between spins
inside a square, and the boundary terms that account for the
interactions between different squares. The contribution of the
latter terms to the free energy vanishes in the limit $L \to \infty$
because the boundary terms grow only linearly in $L$, implying the
existence of the limit defining the free energy.

In our situation, the above argument is not immediately applicable
because we have already taken the scaling limit and are now dealing
with a continuum model. 
We can however try to mimic that argument. For that purpose, we
introduce the functions
\begin{equation} \nonumber
f^t_n(h) \equiv \frac{1}{(2^{n+1})^2} \log{\left(\int \exp(h
\Phi^0_t({\bf 1}_{[-2^n,2^n]^2})) d\nu^t\right)} \, ,
\end{equation}
where $\Phi^0_t$ denotes the near/off-critical magnetization field
with renormalized external field $t$. We now write $\Phi^0_t({\bf
1}_{[-2^n,2^n]^2}) = \sum_k \Phi^0_t(\text{square}_k)$, where
$\text{square}_k$ denotes the $k$th element in a set of equal
squares of fixed size that partition $[-2^n,2^n]^2$. Although the
random variables $\Phi^0_t(\text{square}_k)$ are clearly not
independent, the exponential decay of correlations under $\nu^t$
implies that they are only weakly correlated when the squares are
far apart, suggesting a finite limit for $f^t_n(h)$ as $n \to
\infty$. One can indeed show that the exponential decay of the
covariance between different squares implies that $\limsup_{n \to
\infty} f^t_n(h) < \infty$. The FKG inequality easily implies that
$f^0_n(h) \leq f^t_n(h)$ for $h,t \geq 0$, and that $f^0_n(h)$ and
$f^t_n(h)$ are increasing in $n$. Therefore, one can conclude the
existence of a finite limit for $f^0_n(h)$ as $n \to \infty$.
Comparing the definitions of $f^0_n(h)$ and $f(h)$, this strongly
suggests (and can be used to prove) the existence of the limit
defining $f(h)$.

Integrating~(\ref{eq:conf-cov}), one can check that
\begin{equation} \nonumber
\Phi^0({\bf 1}_{[-\alpha L, \alpha L]^2}) \stackrel{dist.}{=}
\alpha^{15/8} \, \Phi^0({\bf 1}_{[-L,L]^2}) \, ,
\end{equation}
consistent with the scaling law for area measures. If the limit
defining the free energy density exists (and is unique), the above
observation implies that $f(t h)/f(t) = h^{16/15}$, which means that
the free energy density must take the form
\begin{equation} \label{eq:free-energy}
f(h) = C_1 \, h^{16/15}
\end{equation}
for some constant $C_1$.
An immediate consequence of~(\ref{eq:free-energy}) would be the
determination of the tail behavior of the block magnetization:
\begin{equation} \nonumber
Prob(\Phi^0({\bf 1}_{[0,1]^2})>x) \sim \exp{(-C_2 \, x^{16})} \text{
\,\,\,\,\,\,\, for } x>0 \text{ and some constant } C_2>0.
\end{equation}
This result would follow from the methods of~\cite{newman79} (see,
in particular, Theorem 1.4 and Corollary 2.6 there for one-sided
bounds of the same type under similar conditions) and it would show,
incidentally, that the scaling limit magnetization field is not
Gaussian.

\section{Beyond The Ising Model in Two Dimensions} \label{sec:beyond}

In this section, we briefly discuss the applicability of the
approach presented in~\cite{cn-pnas} and in this paper to higher
dimensions, $d>2$, and to $q$-state Potts models with $q>2$.
Although the $d=2$ scaling limit Ising magnetization field $\Phi^0$
should transform covariantly under conformal transformations and
have close connections to the Schramm-Loewner Evolution (SLE), no
conformal machinery seems necessary to establish the existence of
subsequential scaling limits in terms of area measures of critical
FK clusters.

A main ingredient used in this paper is
Prop.~\ref{prop:small-clusters}, which essentially says that
``small'' FK clusters do not contribute to the magnetization in
the scaling limit. This follows from the behavior of the two-point
function at long distance (Prop.~\ref{prop:tau}). Inspecting the
proof, it is easy to check that, in order for
Prop.~\ref{prop:small-clusters} to hold in dimension $d \geq 2$,
$\tau_c(y-x)$ should behave at long distance like
$||y-x||^{-d+2-\eta}$ with $\eta<2$ (see~\cite{cn-pnas}).
Such a decay for $\tau_c$ should be valid for all $d \geq 2$. (In
particular, $\eta$ should be 0 above four dimensions, a result which
has been proved when the number of dimensions is sufficiently
high~\cite{hhs}.) However, there is a significant difference between
dimensions below and above $d=4$, where 4 is the upper-critical
dimension for the Ising model. As we mentioned earlier, for $d=2$
the number of terms in the sum that defines the cutoff
field~(\ref{eq:cutoff-lattice-field}) remains a.s.\ finite in the
scaling limit. This is due to the following result, whose proof is
postponed to the next section.
\begin{proposition} \label{prop:crossprob}
For $z \in {\mathbb R}^2$, let $N^a(z,r_1,r_2)$ denote the number of
distinct clusters ${\cal C}_i^a$ that include sites in both $\{y \in
a \, {\mathbb Z}^2: ||y-z|| < r_1\}$ and $\{y \in a \, {\mathbb
Z}^2: ||y-z|| > r_2\}$.
For any $0<r_1<r_2<\infty$, there exists $\lambda \in (0,1)$ such that
for all $z \in {\mathbb R}^2$ and all small $a>0$ and any $k=1,2,\dots$,
\begin{equation}
\label{eq:crossprob} P_c(N^a(z,r_1,r_2) \geq k) \, \leq \, \lambda^k
\, .
\end{equation}
It follows that for any bounded $D \subset {\mathbb R}^2$ and $\varepsilon>0$,
the number of distinct clusters ${\cal C}_i^a$ of diameter $> \varepsilon$
touching $D$ is bounded in probability as $a \to 0$.
\end{proposition}

The analogue of Prop.~\ref{prop:crossprob} is expected to fail above
the upper-critical dimension $d=4$ (see Appendix A
of~\cite{aizenman1}). When it fails, there can be infinitely many FK
clusters with diameter greater than $\varepsilon$ in a bounded
region and so Prop.~\ref{prop:small-clusters} would not preclude
$\Phi^0$ from being a Gaussian (free) field. But it appears that at
least for $d=3$, both the analogue of Prop.~\ref{prop:crossprob} and
a representation of $\Phi^0_{\varepsilon}$ as a sum of finite
measures with random signs ought to be valid.

An analogous representation for the scaling limit magnetization
fields of $q$-state Potts models also ought to be valid, at least
for values of $q$ such that for a given $d$, the phase transition at
$\beta_c$ is second order. (This was pointed out to the authors
of~\cite{cn-pnas} by J.~Cardy.) The phase transition is believed to
be first order for integer $q \geq 3$ when $d \geq 3$ and for $q >
4$ when $d = 2$ (see~\cite{wu-review}); this leaves, besides the
Ising case, $d = 2$ and $q = 3$ and~$4$. We denote the states or
colors of the $q$-state Potts model by $1,2,\dots,q$, and recall
that in the FK representation on the lattice, all sites in an FK
cluster have the same color while the different clusters are colored
independently with each color equally likely. In the scaling limit,
there would be finite measures $\{\mu_j^{0,q}\}$, and the
magnetization field in the color-$k$ direction would be $\sum_j
\eta_j^k \mu_j^{0,q}$ with the $\eta_j^k$'s taking the value $+1$
with probability $1/q$ (for the color $k$) and the value $-1/(q-1)$
with probability $(q-1)/q$ (for any other color). For a fixed $k$
the $\eta_j^k$'s would be independent as $j$ varies, but for a fixed
$j$ they would be {\it dependent\/} as $k$ varies because $\sum_k
\eta_j^k =0$.

\section{Proofs} \label{sec:proofs}
The proofs of Prop.~\ref{prop:crossprob} and
Prop.~\ref{prop:small-clusters} below follow~\cite{cn-pnas}; we
include them here for completeness. \\

\noindent{\bf Proof of Prop.~\ref{prop:crossprob}.}\,
We define a dual FK model by inserting a bond in the dual lattice,
$({\mathbb Z}^2)^*$, whenever the corresponding dual edge is not
crossed by a bond of the FK configuration on the original lattice,
${\mathbb Z}^2$.

The proof is by induction on $k$. For $k=1$, the result follows from
RSW-type bounds (Theorem~1 of~\cite{dhn}---see~\cite{russo,sewe} for
the original RSW) since $N^a(z,r_1,r_2) \geq 1$ is equivalent to the
{\it absence\/} of a circuit of dual FK bonds (i.e., bonds of the
dual FK model) in the $(r_1,r_2)$-annulus about $z$.
By self-duality at the critical point, this event has the same
probability as the absence of a circuit of FK bonds in the original
FK model, which in turn is bounded away from one as $a \to 0$, by
RSW. Now suppose $N^a(z,r_1,r_2) \geq k-1$. Then one may do an
exploration of the ${\cal C}_i^a$'s that touch $\{y \in a \,
{\mathbb Z}^2: ||y-z|| < r_1\}$ until $k-1$ are found that reach
$\{y \in a \, {\mathbb Z}^2: ||y-z|| > r_2\}$, making sure that all
cluster explorations have been fully completed without obtaining
information about the outside of the clusters.
At that point, the complement $D$ of some random finite $D^c \subset
a \, {\mathbb Z}^2$ remains to be explored and the conditional
random-cluster (FK) distribution in $D$ is $P_c^{{\partial D},F}$
with a {\it free\/} boundary condition on the boundary (or
boundaries) between $D$ and $D^c$. By RSW, the $P_c^{{\partial
D},F}$-probability of a crossing by a sequence of FK bonds in $D$ of
the $(r_1,r_2)$-annulus is bounded above by the original
$P_c(N^a(z,r_1,r_2)\geq 1)$. Thus we have

\begin{eqnarray*}
P_c(N^a(z,r_1,r_2) \geq k) & = & P_c(N^a(z,r_1,r_2) \geq k-1) \\
&   & P_c(N^a(z,r_1,r_2) \geq k|\,N^a(z,r_1,r_2) \geq k-1) \\
& = & P_c(N^a(z,r_1,r_2) \geq k-1) \, E_c[P_c^{{\partial D},F}(N^a(z,r_1,r_2) \geq 1)] \\
& \leq & P_c(N^a(z,r_1,r_2) \geq k-1) \, P_c(N^a(z,r_1,r_2) \geq 1)
\nonumber \\
& \leq &  \lambda ^k \, .
\end{eqnarray*}
The last claim of the proposition follows from~(\ref{eq:crossprob})
because one may choose $O([diam(\Lambda)/\varepsilon]^2)$ points $z_\ell$
in ${\mathbb R}^2$ so that any ${\cal C}_i^a$ of diameter $> \varepsilon$
touching $\Lambda$ will be counted in $N^a(z_\ell, \varepsilon/4,
\varepsilon/2)$ for at least one $z_\ell$. \fbox{} \\

The next proposition corresponds to Hypothesis 1.1 of~\cite{cn-pnas}
(with the exponent $\theta$ there taken to be $1/8$), where it is
shown how, for the critical two-dimensional Ising model, the
hypothesis follows from RSW-type bounds for FK percolation. Such
bounds have recently been proved in~\cite{dhn}. (A derivation of
similar bounds, sufficient to verify Hypothesis 1.1, is also
contained in~\cite{cn-pnas}, but it relies on the convergence of
spin-cluster interfaces to $\text{CLE}_3$, a result that should
follow from Smirnov's work but has not been proved yet.)

\begin{proposition} \label{prop:tau}
There are constants $K_1>0$ and $K_2<\infty$ such that for any small
$\varepsilon >0$ and then for any $x \in {\mathbb Z}^2$ with large
Euclidean norm $||x||$,
\begin{equation} \label{eq:connectivity}
K_2 \tau_c(x_{\varepsilon}) \, \geq \, \tau_c(x) \, \geq \,
K_1 \, \varepsilon^{1/4} \tau_c(x_{\varepsilon}) \,
\end{equation}
for any $x_{\varepsilon} \in {\mathbb Z}^2$ with
$||x_{\varepsilon} - \varepsilon x|| \leq 1/ \sqrt{2}$.
\end{proposition}

\noindent{\bf Proof.}\,
The proposition is an immediate consequence of Prop.~27
of~\cite{dhn}. \fbox{}

\begin{proposition} \label{prop:small-clusters}
For any bounded function $f$ with bounded support,
\begin{equation} \label{eq:small-clusters} \nonumber
\lim_{\varepsilon \to 0} \limsup_{a \to 0}
E_c \left[ \sum_{i:\text{diam}(\gamma^a_i)\leq\varepsilon} (\mu^a_i(f))^2 \right] = 0 \, .
\end{equation}
\end{proposition}

\noindent{\bf Proof.}\,
Using Prop.~\ref{prop:tau}, we can compare $\sum_{z' \in
\Lambda_{\varepsilon' r}} \tau_c(z')$ for small $\varepsilon'$
as $r\to \infty$ to $\sum_{z \in \Lambda_r}\tau_c(z)$ by using
the second inequality of~(\ref{eq:connectivity}) to compare each
$\tau_c(z')$ to the $\tau_c(z)$'s with $\varepsilon' z$ in the
unit length square centered on $z'$ (so that we may take $z'$
as $z_{\varepsilon'}$). Since there are approximately
$(1/{\varepsilon'})^2$ such $z$ sites, we have that
\begin{equation} \label{eq-taucompare} \nonumber
\liminf_{r \to \infty}
\frac{\sum_{z \in \Lambda_r}\tau_c(z)}{({\varepsilon'})^{-7/4}
\sum_{z' \in \Lambda_{\varepsilon' r}} \tau_c(z')}
\, \geq \, K_1 \, .
\end{equation}
Using this lower bound (with $r = 1/2a$ and $\varepsilon' = 2 \varepsilon$)
and~(\ref{eq:Theta}), and letting
$D$ denote the support of $f$ and $D_a \equiv D \cap a \, {\mathbb Z}^2$,
we have that
\begin{eqnarray}
\limsup_{a\to0}
E_c\left[ \sum_{i:\text{diam}(\gamma^a_i)\leq\varepsilon} (\mu^a_i(f))^2 \right]
& \leq & \left(\sup_{x \in D}|f(x)|\right)^2 \limsup_{a\to0} \Theta_a^{2} \,
E_c\left[ \sum_{i:\text{diam}(\gamma^a_i)\leq\varepsilon} |{\cal C}^a_i \cap D|^2 \right]
\nonumber \\
& \leq & \left(\sup_{x \in D}|f(x)|\right)^2 \limsup_{a\to0}
\frac{\sum_{z,w \in D_a, ||z-w|| \leq \varepsilon} \tau_c(w/a-z/a)}{\sum_{x,y \in \Lambda_{1/a}}
\tau_c(y-x)} \nonumber \\
& \leq & \left(\sup_{x \in D}|f(x)|\right)^2
\limsup_{a\to0} \frac{K' (1/a)^2 \sum_{z' \in \Lambda_{\varepsilon /a}}
\tau_c(z')}{K'' (1/a)^2 \sum_{z \in \Lambda_{1/(2a)}} \tau_c(z)} \nonumber \\
& = & K''' \varepsilon^{7/4} \,. \nonumber \,\,\, \fbox{}
\end{eqnarray}

We are now ready to prove the two theorems.\\

\noindent{\bf Proof of Theorem~\ref{thm:subseqential-limits}.}\,
Let $D$ denote the support of $f$; in view of~(\ref{eq:second-moment})
and~(\ref{eq:connectivity}) (compare the proof of
Prop.~\ref{prop:small-clusters}),
\begin{eqnarray*} \nonumber
\limsup_{a \to 0} \left\langle [\Phi^a(f)]^2 \right\rangle_c & = &
\limsup_{a \to 0} E_c\left[ \sum_i (\mu^a_i(f))^2 \right] \\
& \leq & \left(\sup_{x \in D}|f(x)|\right)^2 \limsup_{a \to 0} \Theta_a^2 E_c
\left[\sum_i |{\cal C}^a_i \cap D|^2 \right] < \infty
\end{eqnarray*}
and thus $\Phi^a(f)$ has subsequential limits in distribution as $a \to 0$.
Boundedness of the second moment of $\Phi^a(f)$ and classic Ising
model results (see, e.g., \cite{newman80}~and references therein)
imply that the fourth moment of $\Phi^a(f)$ remains bounded as $a
\to 0$. As a consequence (see, e.g., Problem~14 in Section~8.3
of~\cite{breiman}, p.\ 164), any subsequential limit of $\Phi^a(f)$
has a finite second moment which is the limit of the second moment
of $\Phi^a(f)$. Thus, the distribution of $\Phi^a(f)$ has
subsequential limits in the Wasserstein metric~(\ref{eq:W-distance})
as $a \to 0$.

Since the Euclidean distance makes $[-N,N]^2$ a compact metric
space, the space $C([-N,N]^2)$ of continuous, real-valued functions
on $[-N,N]^2$ with the supremum norm is separable. Every subspace of
a separable metric space is separable, thus the space
$C_0([-N,N]^2)$ of continuous functions with compact support
contained in $[-N,N]^2$ with the supremum norm is also separable.
Any topological space which is the union of a countable number of
separable subspaces is separable, which implies that
$C_0(\mathbb{R}^2) = \bigcup_{N \in {\mathbb N}}C_0([-N,N]^2)$ is
separable. Let $\cal G$ denote a countable, dense subset of
$C_0(\mathbb{R}^2)$; it is clear from the above discussion that we
can choose ${\cal G} = \bigcup_{N \in {\mathbb N}}{\cal G}_N$, where
${\cal G}_N$ is a countable, dense subset of $C_0([-N,N]^2)$. By a
standard diagonalization argument, for every sequence $a_n \to 0$,
there exists a subsequence $a_{n_k} \to 0$ such that, for all $g \in
{\cal G}$, the distribution $P^k_g$ of $\Phi^{a_{n_k}}(g)$ has a
limit $P^0_g \in {\cal P}_2$ in the Wasserstein metric $W_2$ as $k
\to \infty$.




By inspection of the definition of $W_2$, we have the following
straightforward inequalities:
\begin{eqnarray*}
W_2(P^m_f,P^k_f) & \leq & W_2(P^m_f,P^m_g) + W_2(P^m_g,P^k_g) + W_2(P^k_g,P^k_f) \\
& \leq & \left\langle \left|\Phi^{a_{n_m}}(f) - \Phi^{a_{n_m}}(g)\right|^2 \right\rangle_c^{1/2} +
W_2(P^m_g,P^k_g) + \left\langle \left|\Phi^{a_{n_k}}(g) - \Phi^{a_{n_k}}(f)\right|^2 \right\rangle_c^{1/2} \, .
\end{eqnarray*}
Now consider a function $f$ in $C_0({\mathbb R}^2)$ but not in $\cal G$.
Since $f$ has compact support, $f \in C_0([-N_0,N_0]^2)$ for some $N_0$.
If $g \in {\cal G}_{N_0}$, 
the positivity of $\langle S_x S_y \rangle$ for all $x,y$ (or the
independence of the $\eta_i$'s in the FK representation) implies that
\begin{equation*}
\left\langle \left|\Phi^a(f) - \Phi^a(g)\right|^2 \right\rangle_c \leq ||f-g||^2_{\infty} \,
E_c \left[\sum_i \left(\mu^a_i({\bf 1}_{[-N_0,N_0]^2})\right)^2 \right] \, ,
\end{equation*}
and equation~(\ref{eq:second-moment}) and the first inequality
of~(\ref{eq:connectivity}) imply that
$E_c\left[\sum_i \left(\mu^a_i({\bf 1}_{[-N_0,N_0]^2})\right)^2\right]$ is bounded as $a \to 0$.
For $m$ and $k$ sufficiently large, this leads to
\begin{eqnarray*} 
W_2(P^m_f,P^k_f) & \leq & W_2(P^m_g,P^0_g) + W_2(P^0_g,P^k_g) \\
& + & 3 \, ||f-g||_{\infty} \limsup_{a \to 0}
\left( E_c \left[\sum_i \left(\mu^a_i({\bf 1}_{[-N_0,N_0]^2})\right)^2 \right] \right)^{1/2} \, .
\end{eqnarray*}
\noindent (The 3 in the last term is arbitrary, any number greater that 2 would do,
provided that $m$ and $k$ are sufficiently large.)

If $g \in {\cal G}_{N_0}$, as $\ell \to \infty$,
$P^{\ell}_g$ converges to $P^0_g$ in the Wasserstein metric $W_2$
and so the right hand side of the above upper bound for
$W_2(P^m_f,P^k_f)$ can be made arbitrarily small by first choosing
$g$ appropriately, and then taking $m$ and $k$ sufficiently large.
This shows that $P^k_f$ is a Cauchy sequence in $({\cal P}_2,W_2)$.
Since $({\cal P}_2,W_2)$
is complete, as $k \to \infty$,
$P^k_f$ converges in the Wasserstein metric $W_2$ to a probability
distribution $P^0_f \in {\cal P}_2$.

The continuity of
$P^0_{\cdot}:(C_0({\mathbb R}^2),||\cdot||_{\infty}) \longrightarrow ({\cal P}_2,W_2)$
is a consequence of the following inequalities, valid for every $k$,
\begin{eqnarray*} \label{eq:continuity}
W_2(P^0_f,P^0_g) & \leq & W_2(P^0_f,P^k_f) + W_2(P^k_f,P^k_g) + W_2(P^k_g,P^0_g) \\
& \leq & W_2(P^0_f,P^k_f) + \left\langle \left| \Phi^{a_{n_k}}(f) - \Phi^{a_{n_k}}(g) \right|^2 \right\rangle_c^{1/2} + W_2(P^k_g,P^0_g) \\
& \leq & W_2(P^0_f,P^k_f) + ||f-g||_{\infty} \,
\left\langle \left[ \Phi^{a_{n_k}}({\bf 1}_{[-N_0,N_0]^2}) \right]^2 \right\rangle_c^{1/2}
+ W_2(P^k_g,P^0_g) \, ,
\end{eqnarray*}
where $N_0$ is chosen so large that $f,g \in C_0([N_0,N_0]^2)$.
This implies
\begin{equation*}
W_2(P^0_f,P^0_g) \leq ||f-g||_{\infty} \limsup_{a \to 0}
\left\langle \left[ \Phi^a({\bf 1}_{[-N_0,N_0]^2}) \right]^2 \right\rangle_c^{1/2} \,
\end{equation*}
and the conclusion.

We now prove the last statement of the theorem. Since $\cal D$ is a
nuclear space, we can apply the Bochner-Minlos theorem (see for
example~\cite{gj81}, Theorem~3.4.2, p.~52---a proof can be found
in~\cite{gv64}). In order to do so, we define
\begin{equation} \label{eq-charac-func} \nonumber
\chi(f) \equiv \int e^{i x} dP^0_f(x)
\end{equation}
and check the following conditions (where 0 here denotes both
the number 0 and the 0 element of $\cal D$):
\begin{enumerate}
\item Normalization: $\chi(0)=1$,
\item Positivity: $\sum_{k,\ell=1}^m c_k \overline{c_{\ell}} \chi(f_k-f_{\ell}) \geq 0$
for every $m \in {\mathbb N}$, $f_1,\ldots,f_m \in {\cal D}$ and
$c_1,\ldots,c_m \in {\mathbb C}$,
\item Continuity: $\chi(f) \to 1$ as $f \to 0$ (in the topology of $\cal D$).
\end{enumerate}
The first condition is clear from the definition of $\chi$ since
$P^0_f$ is concentrated at the point $x=0$ when $f=0$. To establish
the second condition, let
$F_n \equiv \sum_{k=1}^m c_k e^{i\Phi^{a_n}(f_k)}$
and note that
\begin{equation}
0 \leq \left\langle |F_n|^2 \right\rangle_c = \left\langle
\sum_{k,\ell=1}^m c_k \overline{c_{\ell}} e^{i
\Phi^{a_n}(f_k-f_{\ell})} \right\rangle_c \, . \nonumber
\end{equation}
Along a converging subsequence, $\langle e^{i
\Phi^{a_n}(f_k-f_{\ell})} \rangle_c$ converges to
$\chi(f_k-f_{\ell})$, yielding the desired inequality,
$\sum_{k,\ell=1}^m c_k \overline{c_{\ell}} \chi(f_k-f_{\ell}) \geq
0$.

The remaining step is to establish the continuity of $\chi$. First
note that convergence in the topology of $\cal D$ implies uniform
convergence. With this in mind, the continuity of $\chi$ follows
immediately from the continuity of $P^0_{\cdot}$ proved earlier,
which in particular implies that, if $f$ converges to $g$ uniformly,
the characteristic function of $P^0_f$ converges pointwise to that
of $P^0_g$, and so $\chi(f)$ converges to $\chi(g)$.

In conclusion, by an application of the Bochner-Minlos theorem,
there exists a random, continuous, linear functional $\Phi^0 \in
{\cal D}'$
with characteristic function $\chi$. \fbox{} \\

\noindent{\bf Proof of Theorem~\ref{thm:cutoff-removal}.}\, We first
note that the proof of Theorem~\ref{thm:subseqential-limits} works
also with $\Phi^a(f)$ replaced by $\Phi^a_{\varepsilon}(f)$ for any
$\varepsilon>0$, implying in particular convergence of the
$\varepsilon$-cutoff field in the Wasserstein metric along
subsequences of $a \to 0$. This, combined with a standard
diagonalization argument, implies that for any sequence $a_n \to 0$,
there exists a subsequence $a_{n_k} \to 0$ such that the
distributions of $\Phi^{a_{n_k}}(f)$ and $\Phi^{a_{n_k}}_{1/m}(f)$
converge in the Wasserstein metric $W_2$ as $k \to \infty$ for all
$f \in C_0({\mathbb R}^2)$ and all $m \in {\mathbb N}$. Let $P^0_f$
and $P^0_{f,m}$ denote the respective limits, and let $P^k_f$ denote
the distribution of $\Phi^{a_{n_k}}(f)$ and $P^k_{f,m}$ the
distribution of $\Phi^{a_{n_k}}_{1/m}(f)$.

By inspection of the definition of $W_2$ and the positivity of
$\langle S_x S_y \rangle$ for all $x,y$ (or the independence of
the $\eta_i$'s in the FK representation), we have the following
inequalities:
\begin{eqnarray*}
W_2(P^0_f,P^0_{f,m}) & \leq & W_2(P^0_f,P^k_f) + W_2(P^k_f,P^k_{f,m})
+ W_2(P^k_{f,m},P^0_{f,m}) \\
& \leq & W_2(P^0_f,P^k_f) +
\left\langle \left|\Phi^{a_{n_k}}(f) - \Phi^{a_{n_k}}_{1/m}(f)\right|^2 \right\rangle_c^{1/2}
+ W_2(P^k_{f,m},P^0_{f,m}) \\
& \leq & W_2(P^0_f,P^k_f) +
\left( E_c \left[\sum_{i:diam(\gamma^{a_{n_k}}_i) \leq 1/m} \left(\mu^{a_{n_k}}_i(f)\right)^2 \right] \right)^{1/2}
+ W_2(P^k_{f,m},P^0_{f,m}) \, .
\end{eqnarray*}

The proof of the theorem is concluded by letting first $k \to
\infty$ and then $m \to \infty$, and using the convergence of
$P^k_f$ to $P^0_f$ and of $P^k_{f,m}$ to $P^0_{f,m}$ in the
Wasserstein metric $W_2$, as well as Prop.~\ref{prop:small-clusters}. \fbox{} \\

\noindent{\bf Acknowledgments.}\, The author thanks Christophe
Garban and Charles M. Newman for their suggestions and Ellen Saada
for her patient encouragement during the preparation of this paper.
He also thanks C.M. Newman for many useful comments and discussions,
Wouter Kager for providing Fig.~\ref{fig:loops}, and an anonymous
referee for a careful reading of the manuscript and several useful
suggestions.

\end{document}